
\input amstex
\documentstyle{amsppt}
\NoBlackBoxes
\mag=1200
\pagewidth{5in}
\pageheight{7.3in}
\hcorrection{2mm}
\nologo

\topmatter
\title
Pseudo-Anosov maps and simple closed curves on surfaces
\endtitle
\author
Shicheng Wang$^1$, Ying-Qing Wu and Qing Zhou$^1$
\endauthor
\address
Peking University, Beijing 100871, China
\endaddress
\email 
swang\@sxx0.math.pku.edu.cn
\endemail
\address
University of Iowa, Iowa City, IA 52242, USA
\endaddress
\email
wu\@math.uiowa.edu
\endemail
\address
East China Normal University, Shanghei 200062, China
\endaddress
\email
qzhou\@math.ecnu.edu.cn
\endemail
\thanks  $^1$ Supported by NSF of China.
\endthanks 
\leftheadtext{S-C.~Wang, Y-Q.~ Wu and Q.~Zhou}
\rightheadtext{Pseudo-Anosov maps and simple closed curves}

\abstract 
Suppose $\Cal C$ and $\Cal C'$ are two sets of simple closed curves on
a hyperbolic surface $F$.  We will give necessary and sufficient
conditions for the existence of a pseudo-Anosov map $g$ such that
$g(\Cal C) \cong \Cal C'$.  
\endabstract

\endtopmatter

\document
\define\proof{\demo{Proof}}
\define\endproof{\qed \enddemo}
\define\Int{\text{\rm Int}}
\define\bdd{\partial}
\baselineskip 15pt

Pseudo-Anosov maps are the most important class among surface
homeomorphisms [Th3], and they play important roles in 3-manifold
theory [Th2].

This note is to address the following question: Suppose $F$ is an
orientable closed surface with $\chi(F)<0$, and suppose $c$ and $c'$
are two essential circles on $F$ such that
\roster
\item there is an orientation preserving homeomorphism $f: F \to F$
such that $f(c) =c'$, and
\item $c$ and $c'$ are not in the same isotopy class.  
\endroster 
Does there exists an orientation preserving pseudo-Anosov map $g: F\to
F$ with $g(c)\cong c'$?

Clearly the conditions (1) and (2) are necessary for a positive answer
of the question.  Moreover if we replace ``$F$'' by ``the torus'', and
``pseudo-Anosov'' by ``Anosov'' in the question, then it is known that
the answer is yes. Note also if both $c$ and $c'$ in the
Question are non-separating,
then there is an orientation preserving homeomorphism $f$ such that $f(c)=c'$.

The question is raised by M.~Boileau and S.~Wang to understand degree
one maps between hyperbolic 3-manifolds which are surface bundles over
the circle [BW].  Corollary 2 below gives an affirmative answer to
this question.  We need some definitions in order to state and prove
our main theorem.

A circle $c$ on a compact surface $F$ is {\it essential\/} if $c$ is
not contractible or boundary parallel.  A set $\Cal C= \{c_1, ... ,
c_n\}$ of mutually disjoint circles on $F$ is an {\it independent\/}
set if the curves in $\Cal C$ are essential and mutually non parallel.
Write $\Cal C \cong \Cal C'$ if $\Cal C$ is isotopic to $\Cal C'$.  An
{\it $f$-orbit\/} in $\Cal C$ is a nonempty subset $\Cal C_1$ of $\Cal
C$ such that $f(\Cal C_1) \cong \Cal C_1$.  Notice that if no curve of
$\Cal C$ is isotopic to a curve in $f(\Cal C)$, then $\Cal C$ contains
no $f$-orbit.

\proclaim{Theorem 1} Let $f: F\to F$ be an orientation preserving
homeomorphism on a surface $F$ with $\chi(F) < 0$, and let $\Cal C$ be
an independent set of circles on $F$.  Then there is an orientation
preserving pseudo Anosov map $g: F \to F$ with $g(\Cal C) \cong f(\Cal
C)$ if and only if $\Cal C$ contains no $f$-orbit.  \endproclaim

\proclaim{Corollary 2} Let $c$ and $c'$ be non isotopic essential
curves on a hyperbolic surface $F$ such that $(F, c)$ is orientation
preserving homeomorphic to $(F, c')$.  Then there is an orientation
preserving pseudo Anosov map $g: F \to F$ such that $g(c) \cong c'$.
\endproclaim

\proof This follows immediately from Theorem 1, noticing that if $\Cal
C$ has only one curve $c$ then it contains an $f$-orbit if and only if
$f(c) \cong c$.  \endproof

It is not difficult to see that the condition is necessary.
Below we first assume that $\Cal C$ contains no $f$-orbit, and proceed to
show that there is a pseudo-Anosov map $g: F \to F$ with $g(\Cal C)
\cong f(\Cal C)$.  We remark that it is possible that $\Cal C$
contains no $f$-orbit, and yet $f^k(c) \cong c$ for some $k > 1$, in
which case we can still find a pseudo-Anosov map $g$ such that $g(\Cal
C) \cong f(\Cal C)$.

Given a set of curves $\Cal C$ on $F$, denote by $N(\Cal C)$ a regular
neighborhood of $\Cal C$ in $F$.  The set $\Cal C$ is a {\it
maximal\/} independent set if $\Cal C$ is independent and no component
of $F- \Int N(\Cal C)$ contains essential circles, that is, each
component of $F - \Int N(\Cal C)$ is a pair of pants.  The following
lemma allows us to replace $\Cal C$ in the theorem with a maximal
independent set.

Denote by $\tau_c$ a right hand Dehn twist along a circle $c$ on
$F$.

\proclaim{Lemma 3} Let $f: F\to F$ be a homeomorphism and $\Cal C$ an
independent set of circles on $F$.  If $\Cal C$ contains no $f$-orbit,
then there is a homeomorphism $f' : F\to F$ and a maximal independent
set $\Cal C'$ of circles such that

(1) $\Cal C \subset \Cal C'$,

(2) $f'(c)=f(c)$ for each $c\in \Cal C$, and 

(3) $\Cal C'$ contains no $f'$-orbit.
\endproclaim

\proof Assume that $\Cal C$ is not maximal, so there is a component
$S$ of $F- \Int N(\Cal C)$ with $\chi(S)<0$ which is not a disc with
two holes. Now one can find a pair of essential circles $c'$ and $c''$
in $S$ such that their minimum geometric intersection number is
positive.  Let $\Cal C'= \Cal C\cup c'$.  Then $\Cal C'$ is
independent.  Let $f_{k}= f\circ \tau_{c''}^k$, where $k$ is an
integer. Then for any $k$, $f_{k}(c)=f(c)$ for $c\in C$.  Since the
minimum geometric intersection number of $c'$ and $c''$ is positive,
there are infinitely many isotopy classes of essential circles on $F$
in the family $\{\tau_{c''}^k(c') \, \, | \, \, k\in {\Bbb Z}\}$.
Since there are only finitely many curves in the set $f^{-1}(\Cal
C')$, we may pick $k$ large enough so that $\tau^k_{c''}(c')$ is not
isotopic to any curve in $f^{-1}(\Cal C')$.  Let $f'=f_{k}$ for this
$k$.  Then $f'(c') = f(\tau^k_{c''}(c'))$ is not isotopic to any curve
in $\Cal C'$.

If $\Cal C'$ would contain an $f'$-orbit $\Cal C_1$, then $f'(\Cal
C_1) \cong \Cal C_1$, so $c'\notin \Cal C_1$ because $f'(c')$ is not
isotopic to any curve in $\Cal C'$.  But this would imply that $\Cal
C_1 \subset \Cal C$, and since $f = f'$ on $\Cal C$, it would
contradict our assumption that $\Cal C$ contains no $f$-orbit.  Hence
$\Cal C'$ contains no $f'$-orbit.  The proof now follows by induction
because the number of curves in an independent set is bound above by
$3g + h -3$, where $g$ is the genus and $h$ the number of boundary
components of $F$.  \endproof

Because of Lemma 3, we may assume from now on that the independent set
$\Cal C$ in Theorem 1 is maximal.  Let $M=F\times [0,1]/f$ be the
surface bundle over $S^1$ with fiber $F$ and gluing map $f$, that is,
it is the quotient of $F\times [0,1]$ obtained by identifying $(x,0)$ with
$(f(x),1))$.  Let $q: F\times [0,1]\to M$ be the quotient map.  Denote
by $F^*$ the surface $q(F\times 0)=q(F\times 1)$, by $\Cal C^*$ the
curves $q(\Cal C \times 0)=q(f(\Cal C) \times 1)$ in $M$, and by
$c_i^*$ the curves $q(c_i\times 0) = q(f(c_i)\times 1)$ in $\Cal C^*$,
$i=1,...,n$.  Let $N(\Cal C^*)$ be a tubular neighborhood of $\Cal
C^*$ in $M$; let $M^*= M- \Int N(\Cal C^*)$; let $T_i$ be the torus
$\bdd N(c^*_i)$ on $\bdd M^*$.

Pick a meridian-longitude pair for each $T_i$, with longitude an
intersection curve of $F^* \cap T_i$.  Thus the slopes (i.e isotopy
classes of simple closed curves) on $T_i$ are in one to one
correspondence with the numbers in $\Bbb Q \cup \{\infty\}$.  (See [R]
for more details.)  Define $M^*(q_1,...,q_n)$ to be the manifold
obtained by $q_i$ Dehn filling on $T_i$, $i=1, \ldots, n$.  

\proclaim{Lemma 4} Let $g(k_1,...,k_n) = f \circ \tau_{c_1}^{k_1}\circ
... \circ\tau_{c_n}^{k_n}$. Then

(1) $g(k_1,...,k_n)(c_i)=f(c_i)$ for all $i=1,...,n$ and $k_j\in {\Bbb
Z}$.

(2) $F\times [0,1]/g(k_1,...,k_n)=M^*(1/k_1,...,1/k_n)$.
\endproclaim

\proof (1) Since the curves in $\Cal C$ are mutually disjoint,
$\tau_{c_j} (c_i) = c_i$ for all $i, j$; hence $g(k_1,...,k_n)(c_i) =
f \circ \tau_{c_1}^{k_1}\circ ... \circ\tau_{c_n}^{k_n}(c_i) = f(c_i).$

(2) This is a simple geometric observation: Twisting along the curves
will not change $M^*$, the exterior of the link $\Cal C^*$ in $M$,
while a meridian of $c^*_i$ in $F\times [0,1]/f$ is changed to a
$1/k_i$ curve in $F\times [0,1]/g(k_1, ..., k_n)$.  \endproof

We would like to show that $g(k_1, ... , k_n)$ are isotopic to
pseudo-Anosov maps for some $k_i$.  By a deep theorem of Thurston (see
Otal [Ot] for proof), it suffices to show that $M^*(1/k_1,...,1/k_n)$
is hyperbolic.  It is attempting to apply the hyperbolic surgery
theorem of Thurston [Th1], but we have to be careful because the
manifold $M^*$ may not be hyperbolic.

By an isotopy we may assume that either $f(c_i) = c_j$ or $f(c_i)$ is
not isotopic to $c_j$.  Define an oriented graph $\Gamma$ associated
to $(\Cal C, f)$ as follows.  Each element $c_i$ in $\Cal C$ is a
vertex of $\Gamma$, still denoted by $c_i$, and there is an oriented
edge $c_i c_j$ in $\Gamma$ if and only if $c_j=f(c_i)$.  A subgraph of
$\Gamma$ is a {\it chain\/} if it is homeomorphic to an interval.
Since there is no $f$-orbit and the curves $f(c_i)$ are mutually non
isotopic, we have

\proclaim{Lemma 5} $\Gamma$ has only finitely many components, each of
which is either a single vertex or a chain.  
\endproclaim

For each component $\gamma$ of $\Gamma$, define $C_\gamma$ in $M$ as
follows: If $\gamma$ is an isolated vertex $c_i$, let $C_\gamma$ be
the knot $c^*_i$.  If $\gamma$ has an oriented edge $c_ic_j$, let
$C_{c_ic_j}$ be the annulus $q(c_j \times [0,1])$ in $M$.  Notice that
since $f(c_i) = c_j$, we have
$$q(c_j \times 1) = q(f(c_i) \times 1) = q(c_i \times 0) = c_i^*. $$
Thus $C_{c_ic_j}$ is an annulus with ``top'' boundary curve $c^*_i$
and ``bottom'' boundary curve $c^*_j$.  Now define $C_\gamma$ to
be the union of $C_{c_ic_j}$ over all edges $c_ic_j$ of $\gamma$.
Then $C_\gamma$ is an embedded annulus in $M$ containing all $\{ c^*_i
\,\, | c_i \in \gamma \}$.  Different components of $\Gamma$
correspond to disjoint $C_\gamma$ in $M$.  Let $C_{\Gamma}$ be the
union of $C_\gamma$ for all $\gamma \subset \Gamma$.

\proclaim{Lemma 6} Every incompressible torus in $M^*$ is isotopic to
a torus in $N(C_\Gamma)$.  \endproclaim

\proof Since $F^*$ is incompressible in the surface bundle $M$, the
surface $F' = F^* \cap M^*$ is also incompressible in $M^*$.  Let $T$
be an incompressible torus in $M^*$.  We assume that $T$ has been
isotoped in $M^*$ to meet $F'$ minimally, so $\beta = T \cap F'$
consists of essential curves on both $T$ and $F'$; in particular, it
cuts $T$ into $\pi_1$ injective annuli in $M^*$.  Note that $\beta \neq
\emptyset$, otherwise $T$ would be an incompressible torus in $F\times
[0,1]$; which is impossible as $\chi(F)<0$.

Since $\Cal C$ is maximal, each component of $F'$ is a disk with two
holes, hence each component of $\beta$ is parallel, on $F^*$, to a
unique curve $c^*_i$ in $\Cal C$.  We may isotope $T$ so that $\beta$
lies in $N(C_{\Gamma})$.

An annulus $A$ in $F\times [0,1]$ is called a {\it vertical\/} annulus if
it is isotopic to a product $s\times [0,1]$ for a simple closed curve $s$
in $F$; it is a {\it horizontal\/} annulus if it is boundary parallel.
Each $\pi_1$-injective annulus in $F\times [0,1]$ is either horizontal or
vertical.

Cutting $M^*$ along $F'$, we get a manifold $M'$ homeomorphic to $F\times
[0,1]$.  If a component $A_k$ of $T \cap M'$ is horizontal then its
two boundary components are isotopic to the same $c^*_i$ on $F^*$.
Since $\bdd A_k \subset \beta$ is already in $N(C_\Gamma)$, 
by an isotopy rel $\bdd A_k$ we may push $A_k$ into $N(c^*_i) \subset
N(C_\Gamma)$.  If it is vertical, let $c^*_i = c_i \times 0$ and
$c^*_j = f(c_j) \times 1$ be the two curves on $F\times \bdd [0,1]$
isotopic to $\bdd A_k$.  Then the annulus $A_k$ gives rise to an
isotopy between $c_i$ and $f(c_j)$ on $F$, so by our assumption
above, we have $c_i = f(c_j)$, so $C_{c_ic_j} \subset C_\Gamma$.  It
follows that $A_j$ is also rel $\bdd A_j$ isotopic into $N(C_\Gamma)$.
This completes the proof of Lemma 6.
\endproof

\proclaim{Lemma 7} $X = M - \Int N(C_\Gamma)$ is a hyperbolic
manifold.  \endproclaim

\proof $X$ is irreducible: A reducing sphere $S$ would bound a ball
$B$ in $M$ because $M$, as an $F$ bundle over $S^1$, is irreducible.
Hence $B$ contains some component of $C_\Gamma$.  But since each
component of $C_\Gamma$ contains some essential curve $c_i^*$ of $F$,
and since $F$ is $\pi_1$ injective in $M$, this is impossible.

$X$ is not a Seifert fiber space: Each component of $N(C_\Gamma)$ can
be shrunk to some $N(c^*_i)$, so $X$ contains a nonseparating,
hyperbolic, closed incompressible surface, i.e., the image of $F\times
\frac 12$ under the reverse isotopy of the above shrinking process.
No such surface exists in a Seifert fiber space with boundary because
an essential surface in such a manifold is either horizontal (hence
bounded) or vertical (hence a torus).

$X$ is also atoroidal: If $T$ is an essential torus in $X$, then by
Lemma 6 there is a torus $T'$ in some $N(C_{\gamma}) \cap M^*$ such
that $T \cup T'$ bounds a product region $W = T \times [0,1]$ in
$M^*$.  Since $T$ is essential, $T'$ is $\pi_1$ injective in
$N(C_\gamma) \cap M^* = P \times S^1$, where $P$ is a planar surface.
Note that a $\pi_1$ injective torus in $P \times S^1$ is isotopic to
$c \times S^1$ for some circle $c$ in $P$.  Since $W \subset M^*$, it
contains no component of $\bdd M^*$, so $c$ must be parallel to the
boundary curve of $P$ on $\bdd N(C_\gamma)$, which means that $T'$,
hence $T$, is parallel to $\bdd N(C_\gamma)$.

Since $\partial M$ is a union of tori, the lemma now follows from the
Geometrization Theorem of Thurston for Haken manifolds [Th2].
\endproof

\demo{Proof of Theorem 1}
 If $\Cal C$ contains an $f$-orbit $\Cal C_1$,
then there is a circle $c\in \Cal C$ such that $f^i(c)$ is isotopic to
a curve in $\Cal C_1$ for all $i$, and $f^k(c)\cong c$ for some $k \neq
0$.  Suppose $g: F\to F$ is a homeomorphism such that $g(c') \cong
f(c')$ for all $c' \in \Cal C$, then by induction we have $$g^i(c) =
g(g^{i-1}(c)) \cong g(f^{i-1}(c)) \cong f^i(c)$$ for all $i$.  In
particular, $g^k(c) \cong c$, so $g$ cannot be a pseudo-Anosov map.

Now suppose $\Cal C$ contains no $f$-orbit.  Let $\gamma_1, ...,
\gamma_m$ be the distinct components of $\Gamma$.  After relabeling,
we may assume that $c_i$ is a vertex of $\gamma_i$ for $i \leq m$.
Since $N(C_{\gamma_i})$ is isotopic to $N(c^*_i)$, by performing
trivial surgery on $c^*_j$ for all $j>m$ we get a 3-manifold
$$X = M^*(\emptyset,...,\emptyset, \infty,...,\infty) = M - \Int
N(C_\Gamma).$$ By Lemma 7, $X$ is a hyperbolic manifold, therefore, by
the Hyperbolic Surgery Theorem of Thurston [Th1], $X(1/k_1, ... ,
1/k_m)$ is hyperbolic for sufficiently large $k_i$.  By Lemma 4(2) we
have
$$ X(\frac 1{k_1}, ... , \frac 1{k_m}) = M^*(\frac 1{k_1},...,\frac
1{k_m}, \infty,...,\infty) = F\times [0,1] / g(k_1, ... k_m, 0, ...,
0).$$ By Lemma 4(1) we have $g(k_1, ..., k_m, 0, ..., 0)(\Cal C) =
f(\Cal C)$.  The theorem now follows from Thurston's theorem that
$F\times [0,1] / g$ is hyperbolic if and only if $g$ is isotopic to a
pseudo Anosov map [Th2, Ot].  \endproof

{\eightpoint
ACKNOWLEDGEMENT.  We would like to thank R.D. Edwards, T. Kobayoshi
and W. Thurston for some helpful conversations.}

\enddocument